\documentclass{article}
\usepackage{amsmath,amsthm,amsopn,amstext,amscd,amsfonts,amssymb,verbatim}
\usepackage{color}
\usepackage{natbib}

\def\r{\mathop{\rm r}\nolimits}
\def\diag{\mathop{\rm diag}\nolimits}
\def\vec{\mathop{\rm vec}\nolimits}

\def \build#1#2#3{\mathrel{\mathop{#1}\limits^{#2}_{#3}}}

\newcommand {\boldgreektext}[1] {\boldmath
             \(#1\)\unboldmath}
\newcommand {\boldgreek}[1]
             {\mbox{\boldgreektext{#1}}
            }
\renewenvironment{abstract}
                 {\vspace{6pt}
                  \begin{center}
                  \begin{minipage}{5in}
                  \centerline{\textbf{Abstract}}
                  \noindent\ignorespaces
                 }
                 {\end{minipage}\end{center}}

\theoremstyle{definition}

\newtheorem{rem}{\textbf{Remark}}[section]

\title{\Large \textbf{Some comments about measures, Jacobians and Moore-Penrose inverse}}
\author{
  \textbf{Jos\'e A. D\'{\i}az-Garc\'{\i}a} \thanks{Corresponding author\newline
    }\\
  {\normalsize Universidad Aut\'onoma de Chihuahua} \\
  {\normalsize Facultad de Zootecnia y Ecolog\'{\i}a} \\
  {\normalsize Perif\'erico Francisco R. Almada Km 1, Zootecnia} \\
  {\normalsize 33820 Chihuahua, Chihuahua, M\'exico}\\
  {\normalsize E-mail: jadiaz@uach.mx}\\
  \textbf{Francisco J. Caro-Lopera}\\
  {\normalsize Departament of Basic Sciences} \\
  {\normalsize Universidad de Medell\'{\i}n} \\
  {\normalsize Medell\'{\i}n, Colombia} \\
  {\normalsize E-mail: fjcaro@udem.edu.co} \\[2ex]
}
\date{}
\begin{document}
\maketitle

\begin{abstract}
Some general problems of Jacobian computations in non-full rank matrices are discussed in this work.  In
particular, the Jacobian of the Moore-Penrose inverse derived via matrix differential calculus is
revisited. Then the  Jacobian in the full rank case is derived under the simple and old theory of the
exterior product.
\end{abstract}

\section{Introduction}\label{sec:1}
The multivariate statistical analysis based on singular random matrices is one of the less studied field
of the matrix variate distribution.  The main reason resides in the fact that most of the singular random
matrix distributions do not exist with respect the Lebesgue measure, see \citet{k:68}. Moreover, the
pursued corresponding densities require the computation of Jacobians based on transformation of singular
matrices, which usually do not exist with respect the addressed measure.  At present, only few works are
available around such specialised topic. The area emerged in the 70's with definitions of Wishart and
Beta matrices for non-full rank. Importance and properties of such distributions were not fulfilled until
applications in time series were proposed in \citet{uh:94}, however the results depended on conjectures
about the corresponding Jacobians of leading transformations. Then, the path for a consistent singular
statistics arrived very late with the proof of such conjectures in \citet{dggj:97}, promoting some works
in that research field. In fact, discussions about certain results were needed (\citet{dg:07}). The list
is completed with works focused on random matrix decompositions, generalised inverse and shape theory,
see \citet{dggg:05,dggj:05,dggj:06,dgm:97}.

In this particular case, the problem is increased when the available literature for Jacobians does not
clarify the measures under considerations, and the user can assume that the Lebesgue measure is universal
for such non-full rank matrix transformations. The best example for this issue attains the Moore-Penrose
inverse.

In the full rank case the Jacobian of the Moore-Penrose inverse does not involve a problem because it
exists respect the Lebesgue measure, see \citet{z:1985}, \citet{ns:96}, and \citet{bo:08}, among others.
The addressed works have used a number of modern and classical techniques involving matrix differential
calculus, matrix algebra and statistics.   For the general case, this is in non-full rank matrix case,
\citet{dggj:05,dggj:06} derived the Jacobian and the corresponding explicit Hausdorff measure, by using
factorisation of measures and an approach due to \citet{j:54}. Then the theory of Jacobian computation
via exterior products was set in a general form under a simple idea which can be applied to several
situations.

In this note the Jacobian of the Moore-Penrose inverse is revisited in the general version from matrix
differential calculus. The revisions and discussions are promoted by the standard approach proposed in
\citet{mn:07}. Finally, we complete the exposition by using the James exterior product approach  for
derivation of the Jacobian in the full rank case.

\section{Notation and preliminary results}\label{sec:2}
Let ${\mathcal L}_{m,n}^{+}(q)$ be the linear space of all $n \times m$ real matrices of rank $q \leq
\min(n,m)$ with $q$ distinct singular values. The set of matrices $\mathbf{H}_{1} \in {\mathcal L}_{m,n}$
such that $\mathbf{H}'_{1}\mathbf{H}_{1} = \mathbf{I}_{m}$ is a manifold denoted ${\mathcal V}_{m,n}$,
called Stiefel manifold. In particular, ${\mathcal V}_{m,m}$ is the group of orthogonal matrices
${\mathcal O}(m)$. The rank of a matrix $\mathbf{A}$ is denoted as $\r(\mathbf{A})$ and $\mathbf{A}'$
denotes the transpose matrix of $\mathbf{A}$.

\noindent Observe that, if $\mathbf{X} \in {\mathcal L}_{m,n}^{+}(q)$, we can write $\mathbf{X}$ as
$$
  \mathbf{X} = \left (
          \begin{array}{cc}
            \build{\mathbf{X}_{11}}{}{q \times q} &  \build{\mathbf{X}_{12}}{}{q \times m-q} \\
            \build{\mathbf{X}_{21}}{}{n-q \times q} &  \build{\mathbf{X}_{22}}{}{n-q \times m-q} \\
          \end{array}
          \right ),
$$
such that $r(\mathbf{X}_{11}) = q$. Here $\mathbf{X}_{22}$ is functionally dependent on
$\mathbf{X}_{11}$, $\mathbf{X}_{12}$, $\mathbf{X}_{21}$ by the relation \citep[Problem 1.39, p.54]{g:76}
\begin{equation}\label{X22}
    \mathbf{X}_{22} = \mathbf{X}_{21}\mathbf{X}_{11}^{-1}\mathbf{X}_{12}.
\end{equation}
This is, we shall have $nq+mq -q^{2}$ functionally independent elements in the matrix $\mathbf{X} \in
{\mathcal L}_{m,n}^{+}(q)$, corresponding to the elements of  $\mathbf{X}_{11}, \mathbf{X}_{12}$ and
$\mathbf{X}_{21}$. Moreover, $\mathbf{X}$ can be expressed as
\begin{equation}\label{XX}
  \mathbf{X} = \left (
                  \begin{array}{c}
                    \mathbf{X}_{11} \\
                    \mathbf{X}_{21}
                  \end{array}
               \right )
               \mathbf{X}_{11}^{-1}
               \left (
                  \mathbf{X}_{11},\mathbf{X}_{12}
               \right ).
\end{equation}
Then, without loss of generality, $(d\mathbf{X})$ shall be defined as the exterior product for the
differentials $dx_{ij}$, such that $x_{ij}$ are functionally independent. Explicitly,
\begin{equation}\label{X}
    (d\mathbf{X}) \equiv (d\mathbf{X}_{11})\wedge(d\mathbf{X}_{12})\wedge(d\mathbf{X}_{21}) =
            \bigwedge_{i=1}^{n}\bigwedge_{j=1}^{q}dx_{ij}
            \bigwedge_{i=1}^{q} \bigwedge_{j = q+1}^{m}dx_{ij}.
\end{equation}

\begin{rem}\label{rem1}
Some authors consider that (\ref{X}) define the Lebesgue measure on ${\mathcal L}_{m,n}^{+}(q)$. However,
this measure does satisfy one of the basic properties of the Lebesgue measure about  the invariance under
orthogonal transformations \citep[Theorem 12.2, p. 172]{b:86}, i.e., if $\mathbf{Q} \in {\mathcal O}(m)$
and $\mathbf{H} \in {\mathcal O}(n)$ then, $(d\mathbf{HXQ}) = (d\mathbf{X})$, see \citet{dg:07}.
Therefore, the authors suggest considering a factorisation of a measure on ${\mathcal L}_{m,n}^{+}(q)$.
This idea was developed in \citet{dgm:97}, \citet{dggj:05,dggj:06}, \citet{dggg:05}, and \citet{dg:07}.
\end{rem}

This motivates a revision of some results about the Jacobian and measure with respect can be defined the matrix
transformation $\mathbf{Y} = \mathbf{X}^{+}$, where $\mathbf{X}^{+}$ denotes the Moore-Penrose inverse of
$\mathbf{X}$, see \citet{cm:09}.

When $\mathbf{X}$ has a full rank $\r(\mathbf{X}) = \min(n,m) = r$, the Jacobian was found by
\citet{z:1985}, using algebraic and statistical arguments; and by \citet{ns:96}, with standard matrix
differential calculus. In both cases, $(d\mathbf{X})$ is the Lebesgue measure on ${\mathcal
L}_{m,n}^{+}(r)$, \citet{b:86}. Moreover, without loss of generality, taking $\r(\mathbf{X}) =
m$, hence
\begin{equation}\label{mpicr}
    \mathbf{Y} = \mathbf{X}^{+}= (\mathbf{X}'\mathbf{X})^{-1}\mathbf{X}'
\end{equation}
and they proof that
$$
  (d\mathbf{Y}) = \left|\mathbf{X}'\mathbf{X}\right|^{-n}(d\mathbf{X}),
$$
where $(d\mathbf{X})$ denotes the Lebesgue measure on ${\mathcal L}_{m,n}^{+}(m)$.

In the general case, for non-full rank matrices, i.e., when $\mathbf{X} \in {\mathcal L}_{m,n}^{+}(q)$,
\citet{dggj:05,dggj:06} showed that
$$
  (d\mathbf{Y}) = |\mathbf{D}|^{-2(n+m-q)}(d\mathbf{X}) = \prod_{i = 1}^{q}
    D_{i}^{-2(n+m-q)}(d\mathbf{X}).
$$
Here $\mathbf{X} = \mathbf{H}_{1}\mathbf{D}\mathbf{P}'_{1}$ is the nonsingular part of the
decomposition in singular values of $\mathbf{X}$ \citep{mh:82}, with $\mathbf{D} =
\diag(D_{1},\dots,D_{q})$, $D_{1}> \cdots > D_{q}>0$, $\mathbf{H}_{1} \in {\mathcal V}_{q,n}$ and
$\mathbf{P}_{1} \in {\mathcal V}_{q,m}$. Now $(d\mathbf{X})$ denotes the Hausdorff measure on ${\mathcal
L}_{m,n}^{+}(q)$ and is defined by, see \citet[Section 19]{b:86} and \citet{dgm:97},
\begin{equation}\label{mdx}
    (d\mathbf{X}) = 2^{-q}|\mathbf{D}|^{n+m-2q} \prod_{i<j}^{q}(D_{i}^{2}-D_{j}^{2})(d\mathbf{D})(\mathbf{H}'_{1}
   d\mathbf{H}_{1})(\mathbf{P}'_{1}d\mathbf{P}_{1})
\end{equation}
where $(d\mathbf{D}) = \bigwedge_{i=1}^{q} dD_{ii}$ and $(\mathbf{H}'_{1}d\mathbf{H}_{1})$ defines the
unnormalised invariant probability measure on $V_{q,n}$, see \citet[pp. 67-72]{mh:82}.
Note that the explicit expression for $(d\mathbf{X}) $, provided in (\ref{mdx}) is
not unique. This depends on the factorisation of the measure under consideration. A broad discussion on this
topic can be found at \citet{dggg:05}.

\section{Differentiation and Moore-Penrose inverse}

The context of differentiation techniques can differ strongly from the existence enviroment of the
Jacobian for certain matrix transformation. Statistical text books usually do not clarify the underlying
measures implicit in the computation of a Jacobian.

For example, under certain regularity conditions \citet[Theorem 5, p. 174]{mn:07}, establish that if $\mathbf{X}$ is a
$n \times m$ matrix and $\mathbf{Y} = \mathbf{X}^{+}$, then the matrix of differentials is given by
\begin{eqnarray}
  d\mathbf{Y} &=& -\mathbf{X}^{+}d\mathbf{X}\,\mathbf{X}^{+} + \mathbf{X}^{+}\mathbf{X}^{+'}d\mathbf{X}'
  \left (\mathbf{I}_{n} - \mathbf{XX}^{+}\right) \nonumber\\ \label{mpd}
  && + \left(\mathbf{I}_{m} - \mathbf{X}^{+}\mathbf{X} \right)d\mathbf{X}'\,\mathbf{X}^{+'}\mathbf{X}^{+}.
\end{eqnarray}
where $d\mathbf{A}$ denotes the matrix of differentials of $\mathbf{A}$. This result was proposed
originally by \citet{gp:73} in terms of Fr\'echet derivatives of orthogonal projector associated to
$\mathbf{X}$ and $\mathbf{X}^{+}$.

Now, for revisiting the classical exposition about the existence of certain measures involved in the
standard theory of Jacobian computation,  we follow the approach proposed by \citet{mn:07} in order to
find the Jacobian $\mathbf{Y} = \mathbf{X}^{+}$. If we just follow the technique, then by applying
\citet[Eq. (5), p. 35]{mn:07} and observing that $\mathbf{XX}^{+}$ and $\mathbf{X}^{+}\mathbf{X}$  are
symmetric matrices, we obtain
\begin{eqnarray*}
  d\vec \mathbf{Y} &=& - \left(\mathbf{X}^{+'} \otimes \mathbf{X}^{+} \right) d\vec \mathbf{X} +
  \left[\left (\mathbf{I}_{n} - \mathbf{XX}^{+}\right) \otimes \mathbf{X}^{+}\mathbf{X}^{+'}\right]d\vec
  \mathbf{X}'\\
   && + \left [\mathbf{X}^{+'}\mathbf{X}^{+} \otimes \left(\mathbf{I}_{m} - \mathbf{X}^{+}\mathbf{X} \right)\right ]d\vec
  \mathbf{X}'.
\end{eqnarray*}
Here $\vec$ denotes the vectorisation operator and $\otimes$ denotes the Kronecker product, see
\citet[Section 2.2 ]{mh:82}. Recalling that $\mathbf{K}_{mn}\vec \mathbf{A} = \vec \mathbf{A}'$, where
$\mathbf{K}_{mn}$ denotes the commutation matrix \citep[Section 7, p. 54]{mn:07}, we get
\begin{eqnarray}
  d\vec \mathbf{Y} &=& \left \{- \mathbf{X}^{+'} \otimes \mathbf{X}^{+} +
  \left[\left (\mathbf{I}_{n} - \mathbf{XX}^{+}\right) \otimes \mathbf{X}^{+}\mathbf{X}^{+'}\right.
  \right. \nonumber\\ \label{dvX}
  && + \left .\left .\left(\mathbf{I}_{m} - \mathbf{X}^{+}\mathbf{X} \right) \otimes \mathbf{X}^{+'}
  \mathbf{X}^{+}\right ]\mathbf{K}_{mn}\right \} d\vec \mathbf{X},
\end{eqnarray}
Then, by applying the first identification theorem for matrix function, see \citet[Eqs. (3) and (4), p.
198]{mn:07} we finally obtain
\begin{eqnarray}
(d\mathbf{Y}) &=& |J(\mathbf{X} \rightarrow \mathbf{Y})|(d\mathbf{X}) \nonumber\\
  &=&\left |- \mathbf{X}^{+'} \otimes \mathbf{X}^{+} +
  \left[\left (\mathbf{I}_{n} - \mathbf{XX}^{+}\right) \otimes \mathbf{X}^{+}\mathbf{X}^{+'}\right.
  \right. \nonumber\\ \label{measure1}
  && + \left . \left .\left(\mathbf{I}_{m} - \mathbf{X}^{+}\mathbf{X} \right) \otimes \mathbf{X}^{+'}
  \mathbf{X}^{+}\right ]\mathbf{K}_{mn}\right| (d\mathbf{X}),
\end{eqnarray}
where $| \cdot |$ denotes the determinant.

After such standard procedure, we ask from (\ref{measure1}): What is the arising measure $(d\mathbf{X})$?
The emerging problem comes from the computation of the Jacobian of the transformation $\mathbf{Y} =
\mathbf{X}^{+}$. Because the Jacobian exists with respect to the Lebesgue measure if the elements of the
matrix $\mathbf{X}$ are functionally independent real variables, see \citet{m:97}. But (\ref{X22})
definitely shows that this hypothesis is not fulfilled. Then, the Jacobian (\ref{measure1}) obtained with
such standard method is incorrect and/or the corresponding measure is not the Lebesgue measure.

Jacobian computation by using linear structure theory lead to the same question if we use certain
analogies for Jacobians of transformations involving functions of symmetric, triangular or diagonal
matrices. For example, calculation of the same Jacobian, via \citet{m:88}, can suggest the following
procedure: propose a matrix, say $\mathbf{M}_{mn}$, such that $\mathbf{M}_{mn} \vec \mathbf{X} =
\vec\mathbf{X}_{I}$, where $\vec(\mathbf{X}_{I})$ denotes the vectorisation of $\mathbf{X}$, but only
considering the functionally independent elements in $\mathbf{X}$. Then we can proceed as in the case of
symmetric, triangular or diagonal matrices. Unfortunately, this is not possible either, since
$\mathbf{X}_{22}$ is not a linear function of the remaining functionally independent elements in
$\mathbf{X}$, see equation (\ref{X22}).

Then, the correct Jacobian of the transformation $\mathbf{Y} = \mathbf{X}^{+}$ and its
corresponding measure, would require a similar expression  to (\ref{dvX}), but in terms only of
functionally independent elements in $\mathbf{X}$. The method shall be sketched in the next few lines. For
The first stage must propose an explicit expression of $\mathbf{X}^{+}$, but only in terms
of the functionally independent elements in $\mathbf{X}$, that is in terms of $\mathbf{X}_{11}$,
$\mathbf{X}_{12}$, and $\mathbf{X}_{21}$. This is a feasible objective because (\ref{XX}) and
\citet[Corollary 1.4.2, p. 22]{cm:09} turn into
\begin{equation}\label{mpigI}
    \mathbf{X}^{+} = \left (
                        \begin{array}{c}
                          \mathbf{X}_{11}^{'} \\
                          \mathbf{X}_{12}^{'}
                        \end{array}
                     \right )
                     \left(
                       \mathbf{X}_{11}\mathbf{X}_{11}^{'}+\mathbf{X}_{12}\mathbf{X}_{12}^{'}
                     \right )^{-1} \mathbf{X}_{11}
                     \left(
                       \mathbf{X}_{11}^{'}\mathbf{X}_{11}+\mathbf{X}_{21}^{'}\mathbf{X}_{21}
                     \right )^{-1}
                     \left(
                       \mathbf{X}_{11}^{'},\mathbf{X}_{21}^{'}
                     \right ).
\end{equation}
And finally, the computation ends by expressing the differentials of (\ref{mpigI}) in terms of
$(d\mathbf{X}_{11})$, $(d\mathbf{X}_{12})$ and $(d\mathbf{X}_{21})$. As first sight this task  seems
cumbersome and laborious; and according to Remark \ref{rem1} the desired result should be difficult to
achieve.

\begin{rem}
The reader should not mislead that (\ref{X}) is not a possible measure, we just simply point out that neither the Lebesgue measure nor the Hausdorff measure are of type  (\ref{X}), because both are invariant under orthogonal transformations, see \citet[Theorem 19.2, p. 252]{b:86}.
\end{rem}

For completeness, in the full rank case,  if $\r(\mathbf{X}) = q = m$, we have that
$\mathbf{X}^{+}\mathbf{X} = \mathbf{I}_{m}$, then
\begin{eqnarray*}
  (d\mathbf{Y}) &=& \left |- \mathbf{X}^{+'} \otimes \mathbf{X}^{+} + \left [\left (\mathbf{I}_{n} -
  \mathbf{XX}^{+}\right) \otimes \mathbf{X}^{+}\mathbf{X}^{+'}\right ]\mathbf{K}_{mn}\right|
  (d\mathbf{X}) \\
   &=& |\mathbf{X}'\mathbf{X}|^{-n}(d\mathbf{X}).
\end{eqnarray*}
In this case $(d\mathbf{X})$ does denote the Lebesgue measure on ${\mathcal L}_{m,n}^{+}(m)$ and  agrees
with \citet[equation (3)]{ns:96}. Similarly, $\r(\mathbf{X}) = q = n$, we have that
$\mathbf{X}\mathbf{X}^{+} = \mathbf{I}_{n}$, then
\begin{eqnarray*}
  (d\mathbf{Y}) &=& \left |- \mathbf{X}^{+'} \otimes \mathbf{X}^{+} +
  \left[\left(\mathbf{I}_{m} - \mathbf{X}^{+}\mathbf{X} \right) \otimes \mathbf{X}^{+'}
  \mathbf{X}^{+}\right ]\mathbf{K}_{mn}\right| (d\mathbf{X}) \\
   &=& |\mathbf{XX}'|^{-m}(d\mathbf{X})
\end{eqnarray*}
where now $(d\mathbf{X})$ denotes the Lebesgue measure on ${\mathcal L}_{m,n}^{+}(n)$.

Finally, instead of using the modern theory of linear structures or algebraical approaches for the full
rank case Jacobian computation given by  \citet{z:1985} and \citet{ns:96},  we propose the simple and
elegant old ideas via exterior products based on James works of the 1950's.

First  define $\mathbf{A} = \mathbf{Y}'\mathbf{Y}$ and $\mathbf{B} =
\mathbf{X}'\mathbf{X}$ in (\ref{mpicr}). Then
\begin{eqnarray*}
  (d\mathbf{Y}) &=& 2^{-m} |\mathbf{A}|^{(n-m-1)/2}(d\mathbf{A})\wedge (\mathbf{H}_{1}^{'}d\mathbf{H}_{1})\\
  (d\mathbf{X}) &=& 2^{-m} |\mathbf{B}|^{(n-m-1)/2}(d\mathbf{B})\wedge (\mathbf{G}_{1}^{'}d\mathbf{G}_{1}).
\end{eqnarray*}
where $\mathbf{H}_{1}, \mathbf{G}_{1} \in \mathcal{V}_{m,n}$. This Jacobian was found by \citet{j:54},
\citet{h:55} and \citet{ro:57} via singular value, polar and QR factorisation, respectively. Then
\begin{eqnarray}
  \label{11}
  (d\mathbf{A}) &=& 2^{m} |\mathbf{A}|^{-(n-m-1)/2}(d\mathbf{Y})\wedge (\mathbf{H}_{1}^{'}d\mathbf{H}_{1})^{-1}\\
  \label{12}
  (d\mathbf{B}) &=& 2^{m} |\mathbf{B}|^{-(n-m-1)/2}(d\mathbf{X})\wedge (\mathbf{G}_{1}^{'}d\mathbf{G}_{1})^{-1}.
\end{eqnarray}
Also, note that
$$
  \mathbf{A} = \mathbf{Y}'\mathbf{Y}= (\mathbf{X}'\mathbf{X})^{-1} \mathbf{X}'\mathbf{X} (\mathbf{X}'\mathbf{X})^{-1}
  = (\mathbf{X}'\mathbf{X})^{-1} = \mathbf{B}^{-1}.
$$
Then, by \citet[Theorem 2.1.8, p.59]{mh:82},
\begin{equation}\label{da}
    (d\mathbf{A})= |\mathbf{B}|^{-(m+1)}(d\mathbf{B}).
\end{equation}
Now, substitute (\ref{12}) in (\ref{da}) and match the result  with (\ref{11}). Then, by the uniqueness
of the nonnormalised measure on Stiefel manifold, $(\mathbf{H}'_{1}d\mathbf{H}_{1}) =
(\mathbf{G}'_{1}d\mathbf{G}_{1})$ \citep[Section 2.1.4]{mh:82}, the required result is obtained.
\begin{eqnarray*}
  (d\mathbf{Y}) &=& |\mathbf{A}|^{(n-m-1)/2} |\mathbf{B}|^{-(m+1)}|\mathbf{B}|^{-(n-m-1)/2} (d\mathbf{X})\\
   &=& |\mathbf{X}'\mathbf{X}|^{-(n-m-1)/2} |\mathbf{X}'\mathbf{X}|^{-(m+1)} |\mathbf{X}'\mathbf{X}|^{-(n-m-1)/2}
   (d\mathbf{X}) \\
   &=& |\mathbf{X}'\mathbf{X}|^{-n} (d\mathbf{X}).
\end{eqnarray*}

\end{document}